\documentclass[a4paper,12pt]{article}
\usepackage[utf8]{inputenc}

\usepackage{amsmath,amsthm}
\usepackage{amssymb}
\usepackage{xcolor}

\newcommand{\R}{\mathbb{R}}
\newcommand{\N}{\mathbb{N}}
\newcommand{\D}{\mathbb{D}}
\newcommand{\de}{\partial}
\newcommand{\ep}{\epsilon}

\newcommand{\linf}{L_\infty(\R_+;L_\infty)}
\renewcommand{\qed}{\begin{flushright} $\square$ \end{flushright}}

\newtheorem{thm}{Theorem}
\newtheorem{lem}{Lemma}

\newtheorem{cor}{Corollary}

\renewcommand{\div}{{\rm div\,}}
\title{Reacting multi-component fluids -- regular solutions in Lorentz spaces}
\author{Piotr Bogusław Mucha$^*$  \& Tomasz Piasecki\thanks{Institute of Applied Mathematics and Mechanics, University of Warsaw, Banacha 2, 02-097 Warsaw, Poland}}

\begin{document}

\maketitle

\centerline{{\em Dedicated to Professor Yoshihiro Shibata on occasion of his 70th anniversary}}

\begin{abstract}
    The paper deals with the analysis of a model of a multi-component fluid 
    admitting chemical reactions. The flow is considered in the incompressible regime. 
    The main result shows the global existence of regular solutions under the assumption of suitable smallness conditions. In order to control the solutions a special structure condition on the derivatives of chemical production functions determining the reactions is required. The existence is shown in a new critical functional framework of Lorentz spaces of type $L_{p,r}(0,T;L_q)$, which 
    allows to control the integral $\int_0^\infty \|\nabla u(t)\|_{\infty} dt$.
\end{abstract}


\section{Introduction}
 
We are interested in the existence of regular solutions to a system describing the flow of a mixture of incompressible, reacting constituents. In the chosen setting the components are given by their fractional densities $\rho_i$ and share a common velocity field. The system consists of the  classical incompressible Navier-Stokes system coupled with a system of nonhomogeneous transport equations describing the evolution of fractional masses and their change under chemical reactions.
The coupling is realized through the dependence of the viscosity coefficient on fractional masses. Precisely, the system reads:
\begin{equation} \label{sys}
\begin{array}{lcl}
 \rho_{i,t} + u \cdot \nabla \rho_i = r_i(\vec \rho\,),\quad i=1,\ldots,M & \mbox{ in } & \R^3 \times [0,T),\\[5pt]
 \rho u_t + \rho u\cdot \nabla u - \div \!\!\left(\nu (\vec  \rho\,) \D(u)\right) + \nabla \pi =0& \mbox{ in } & \R^3 \times [0,T),\\[5pt]
 \div u=0 & \mbox{ in } & \R^3 \times [0,T),
\end{array}
\end{equation}
where
\begin{equation}
 \rho = \sum_{i=1}^M \rho_i, \quad \vec \rho = (\rho_1, \rho_2, ... ,\rho_M).
\end{equation}
In the system we look for the velocity $u$, the vector of fractional densities of components of the fluid $\vec \rho$ and the internal pressure $\pi$. $M$ denotes the number of constituents, $\nu(\cdot)$ denotes the viscosity coefficient which is assumed to be a function of the whole $\vec \rho$. The reactions are coded in the production functions $r_i(\cdot)$.
The system \eqref{sys} is supplied with initial conditions 
\begin{equation} \label{init}
\rho|_{t=0}=\rho^0, \quad 
\vec \rho\,|_{t=0}={{\vec \rho}}_0, \quad u|_{t=0}=u_0 
\mbox{ \ \ with \ \ } 
\vec \rho_0 = (\rho_1^0,\ldots,\rho_m^0) \quad 
\rho^0 = \sum_{i=1}^M \rho_i^0.
\end{equation}

The analysis of the properties of flows of mixtures has been attracting the attention of the community of mathematical fluid mechanics in recent years. 
A widely accepted mathematical model of a mixture is the Maxwell-Stefan system
which describes diffusion of constituents. This system is usually degenerate parabolic, however it reveals an entropy structure which allows to show well-posedness of the problem. A nice self-contained introduction to this approach is given in \cite{Jungel}. Mathematical properties of a pure Maxwell-Stefan reaction-diffusion system has been investigated among others in \cite{B2010},\cite{HMPW13}.

Another approach is to look at the mixture as a fluid. If we assume that the components share the same velocity then such situation is described by Maxwell-Stefan type systems coupled with the Navier-Stokes equations. 
A rigorous derivation of a class of such models has been shown in the monograph \cite{Gio}, where also the existence of regular solutions close to equilibrium is shown in the whole space. 
In case of compressible flow the coupling is usually assumed to be realised in the form of the pressure, which depends on fractional densities. In the incompressible case the pressure is unknown, therefore a natural way is to assume coupling in a variable viscosity coefficient. Such model was introduced in \cite{BP1} and is assumed also in this paper.

The simplest assumption on the diffusion matrix is the Fick Law where the diffusion matrix is diagonal, therefore lack of cross-diffusion is assumed. The existence of weak solutions for such model coupled with compressible Navier-Stokes equations has been shown in \cite{FPT}. 
If we assume cross-diffusion the analysis becomes more delicate \cite{MPZ2}. Global existence of weak solutions under additional assumptions on the relation between the viscosity coefficients has been shown in \cite{MPZ}. 

Strong solutions for a general model with cross-diffusion on domains has been shown to exist locally in time and globally under additional smallness assumptions in \cite{PSZ1} for the case of two constituents. The approach is based on $L_p-L_q$ maximal regularity. The result has been generalized to arbitrary number of constituents in \cite{PSZ2}-\cite{PSZ3}.
Similar results have been shown recently under different structural assumptions on the diffusion subsystem in \cite{BoDr}.

All above results essentially exploit the regularizing effect of diffusion in the system describing the evolution of fractional masses. To the best knowledge of the authors, the issue of well-posedness for a Navier-Stokes system coupled with a transport-reaction system has not been dealt with yet. 

In the present paper we restrict ourselves to incompressible mixture, leaving the compressible case without diffusion in the species subsystem as an interesting open problem for future research. We aim at construction of global in time regular solutions, hence we are required to restrict our analysis to small data case.

Mathematical analysis of incompressible multicomponent flows dates back to a two component model with Fick law in \cite{BdV1} for inviscid fluid and in \cite{BdV2}-\cite{BdV3} in the viscous case.
Global existence of weak solutions to the incompressible Navier-Stokes-Stefan-Maxwell system with arbitrary data was proved   
independently using different techniques in \cite{CJ13} and \cite{MT13}. 
Strong solutions in the $L_p-L_q$ maximal regularity were investigated in already mentioned paper \cite{BP1}. 

\medskip 

{\bf Structural assumptions on the chemical reactions.}
Our goal is to investigate a system describing the flow of a mixture of constituents which undergo chemical reactions and the average velocity and the total density obey the classical Navier-Stokes system, however without diffusion in the subsystem corresponding to balance of fractional masses. It is well known that diffusion gives additional regularity and some decay in time. Without this regularizing effect we have to investigate carefully the properties of the species subsystem to obtain estimates necessary to show the well-posedness of the problem for large time. 

For this reason we need to introduce several assumptions on the structure of chemical reactions. We assume that reactions take place in a dilutant denoted by $w$ and fractional masses of other constituents are small compared to the mass of the latter. Furthermore, we consider a single reaction in which $k$ reactants $a_{1},\ldots,a_{k}$ give $l$ products $b_{1},\ldots,b_{l}$. 
We assume that production of each constituent $b_{j}$ is a known, constant fraction $\theta_j$ of the whole production (which is a natural assumption if we consider some known chemical reaction). An important assumption here is that the chemical production rates $r_k$ depend only on the reactants, and not on the products. This assumption is justified especially that we assume the presence of the dilutant which dominates in the quantitative sense other constituents. Formulating precisely the above assumptions, we have 
\begin{equation} \label{def:vecrho}
\vec \rho = (w,a_1,\ldots,a_k,b_1,\ldots,b_l), \quad M=1+k+l,   
\end{equation}
and denoting 
\begin{equation} \label{def:vecab}
\vec a = (a_1,\ldots,a_k), \qquad 
\vec b = (b_1,\ldots,b_l)
\end{equation}
we  rewrite subsystem \eqref{sys}$_1$ as 
\begin{equation} \label{subsys}
\begin{aligned}
&w_t+u\cdot \nabla w=0, \\[3pt]
&a_{i,t} + u \cdot \nabla a_i = - \omega_i(\vec a), \quad i=1,\ldots,k ,\\[3pt]
&b_{j,t} + u \cdot \nabla b_j = \theta_j \sum_{i=1}^k \omega_i(\vec a), \quad j=1,\ldots,l,
\end{aligned}
\end{equation}
where $\omega_i$ are given nonnegative functions. In order to keep the sign of the densities we require that
\begin{equation} \label{omega:zero}
\omega_i(\vec a)\geq 0, \quad \omega_i(a_1,...,a_{i-1},0,a_{i+1},...,a_k)=0.
\end{equation} 
Furthermore, $\theta_j$ are nonnegative constants satisfying 
$\sum_{j=1}^l \theta_j=1.$  The form of (\ref{def:vecrho})
assumes that we consider just one reaction, the model can be extended to more reactions but we do not want to check the patience of the readers. From the modelling and mathematical viewpoint such change makes the system more complex in description only.

\section{Functional spaces and the main result}

We work in the setting of Lorentz spaces in time following the ideas from a recent work \cite{DMT}, where such functional framework is used in the context of simple compressible dynamics of a viscous flow. We show that similar estimates yield the well posedness of the model of incompressible mixture without the regularizing effect of species diffusion under assumptions on the structure of species subsystem described above and certain assumption on chemical production functions $\omega_i$ which we formulate below in the statement of our main result. We also present a special case of quite simple form of reaction functions fulfilling these requirements. Before formulating our main result we shall introduce the notation used in the paper.

Mostly we work in the Lebesgue spaces 
$L_p(X,\mu)$ setting generating Besov and Lorentz ones 
by interpolation. The latter are not too popular in the analysis of problems in PDEs, although their properties
are sometimes more interesting than the features of Besov ones.

\smallskip 

{\bf Lorentz spaces} can be defined on a measure space $(X,\mu)$ by real interpolation between Lebesgue spaces as 
\begin{equation*}
L_{p,r}(X,\mu) = ( L_\infty(X,\mu),L_1(X,\mu))_{1/p,r} \quad {\rm for} \quad
p\in(1,\infty),\; r\in [1,\infty].
\end{equation*}
In a direct way the spaces are defined as spaces of functions for which the norm 
\begin{equation*}
    \|f\|_{L_{p,r}}=\left\{
    \begin{array}{lr}
   \displaystyle p^{1/r} \left( \int_0^\infty (s|\{|f|>s\}|^{1/p})^r \frac{ds}{s}\right)^{1/r}     &  \mbox{ for \ } r<\infty, \\[14pt]
    \displaystyle \sup_{s>0} s|\{ |f|>s\}|^{1/p}     &  \mbox{ for \ } r=\infty,
    \end{array}
    \right. 
\end{equation*}
where $|\cdot|$ denotes $\mu(\cdot)$, is finite. The
factor $p^{1/r}$ ensures that $\|f\|_{L_{p,p}}=\|f\|_{L_p}$.

Let us recall several properties of Lorentz spaces which can be found in \cite[Chapter 4.4]{BS}, see also \cite{G}.  
For brevity we write $L_p, L_{p,r}$ instead of $L_p(X,\mu), L_{p,r}(X,\mu)$. 
\begin{enumerate}
    \item Imbedding: \begin{equation}\label{imbed}
    L_{p,p}=L_p, \; L_{p,r_1} \subset L_{p,r_2} \quad {\rm for} \quad r_1 \leq r_2.
    \end{equation}
    \item H\"older inequality: for $1<p,p_1,p_2<\infty$ and $1\leq r,r_1,r_2\leq\infty$ such that $\frac{1}{p}=\frac{1}{p_1}+\frac{1}{p_2}$ and $\frac{1}{r}=\frac{1}{r_1}+\frac{1}{r_2}$ we have   
    \begin{equation} \label{holder}
        \|fg\|_{L_{p,r}} \leq C \|f\|_{L_{p_1,r_1}} \|g\|_{L_{p_2,r_2}}.     
    \end{equation}
    Inequality holds also for $(p_1,r_1)=(1,1)$ and $(p_1,r_1)=(\infty,\infty)$, where $L_{1,1}=L_1$ and $L_{\infty,\infty}=L_{\infty}$.
    \item If $f \in L_\infty(\Omega)$ and $g \in L_{p,q}(\Omega)$ then $fg \in L_{p,q}(\Omega)$ for all $p,q \in [1,\infty]$ and 
    \begin{equation} \label{0:1}
        \|fg\|_{L_{p,q}}\leq \|f\|_\infty \|g\|_{L_{p,q}}.
    \end{equation}
    \item For $\alpha>0$, $p,q,p\alpha,q\alpha \in [1,\infty]$
    \begin{equation} \label{0:2}
        \|f^\alpha\|_{L_{p,q}} \sim \|f\|_{L_{p\alpha,q\alpha}}^\alpha.
    \end{equation}
\end{enumerate}

{\bf Besov spaces} are split into two subclasses, homogeneous and nonhomogenenous one, we denote them respectively by $\dot B^s_{p,q}(\Omega)$ and $B^s_{p,q}(\Omega)$. The most compact definition  can be given  by real interpolation of classical Sobolev spaces.  For $s\in (0,1)$ we have
\begin{equation*}
\dot B^s_{p,q}(\Omega)=(L_p(\Omega);\dot W^1_p(\Omega))_{s,q}, \qquad    
B^s_{p,q}(\Omega)=(L_p(\Omega);W^1_p(\Omega))_{s,q}.
\end{equation*}
In general the Besov spaces are the interpolation family, i.e.
\begin{equation} \label{0:3a}
    \dot B^s_{p,q}=(\dot B^{s_1}_{p,q_1},\dot B^{s_2}_{p,q_2})_{\theta,q}, \qquad
     B^s_{p,q}=(B^{s_1}_{p,q_1},B^{s_2}_{p,q_2})_{\theta,q}, 
\end{equation}
with $\theta \in (0,1)$, $s=(1-\theta)s_1+\theta s_2$ and $s_1,s_2 \in \R$, $p,q,q_1,q_2\in [1,\infty]$.

Here we recall that $\dot W^m_p(\R^d)$ is the homogeneous Sobolev space supplemented with a seminorm
\begin{equation*}
    \|f\|_{\dot W^m_p} = \sum_{|\alpha|=m} \|\partial^\alpha f\|_{L_p}.
\end{equation*}
Note that it is a linear space but not necessary a Banach one (for $mp>d$ it is not).

A direct definition of Besov spaces is usually described by the Fourier transform. As we do not use this language, 
we omit it here and refer to \cite{D-book},
where one can find also proofs of basic results for these spaces.
We recall two of them, the first one is the imbedding theorem
\begin{equation} \label{0:3}
    \dot B^s_{p,r} \subset \dot B^{s-d(\frac{1}{p}-\frac{1}{q})}_{q,r}, \mbox{ in particular \ }
    \dot B^{d/p}_{p,1} \subset \dot B^0_{\infty,1} \subset L_\infty \mbox{ \ \ for \ } p,q,r \in [1,\infty],
\end{equation}
where $d$ is the space dimension.
The second is the interpolation inequality: for any $s\in \R$, $\theta \in (0,1)$, $s=(1-\theta)s_1+\theta s_2$ and $s_1,s_2 \in \R$, $p,q,q_1,q_2\in [1,\infty]$
\begin{equation} \label{0:4}
    \|f\|_{\dot B^{s}_{p,q}} \leq C\|f\|_{\dot B^{s_1}_{p,q_1}}^{1-\theta} \|f\|_{\dot B^{s_2}_{p,q_2}}^\theta.
\end{equation}


\medskip

Let us also introduce a brief notation for a regularity class in which we shall work: 
\begin{equation} \label{def:3}
W^{2,1}_{p,(q,r)}(\R^3 \times \R_+):= \{ z \in C_b(\R_+;\dot B^{2-2/q}_{p,r}(\R^3)): \; \de_t z,\nabla^2 z \in L_{q,r}(\R_+;L_p(\R^3) \}
\end{equation}
with the norm
\begin{equation*}
    \|z\|_{W^{2,1}_{p,(q,r)}} :=
    \|z\|_{L_\infty(\R_+;\dot B^{2-2/q}_{p,r}(\R^3))} +
    \|\partial_t z, \nabla^2 z\|_{L_{q,r}(\R_+;L_p(\R^3))}.
\end{equation*}
Note that the above space is homogeneous, although we do not put the dot over $W$. 
Next, by $\phi(\cdot)$ we shall denote a continuous function (possibly of many variables) such that $\phi(0)=0$ (therefore it plays a role of a small constant dependent on the data). 
We also sometimes skip the spatial domain in notation of norms, for instance we write $L_p$ instead of $L_p(\R^3)$. Since the space domain is always $\R^3$ in our considerations, it should not lead to misunderstanding. Finally we will write $A \lesssim B$ if $A \leq CB$ where $C$ is a generic constant.

\medskip
\noindent
We are now in a position to formulate our main result
\begin{thm} \label{t1}
Assume the species subsystem \eqref{sys}$_1$ has the form \eqref{subsys} with chemical production rates $\omega_j \in C^1(\R^k)$ satisfying \eqref{omega:zero} and the following structural assumption for $p=3$ and $6$
\begin{equation}
%
\sum_{m=1}^k a_m^{-(p-1)\alpha_p}|a_{m,x_j}|^{p-2}a_{m,x_j} \de_{x_j} \big( \omega_m(\vec a)a_m^{-\alpha_p} \big) \geq 0 
 \label{omega:pos}
\end{equation}
for $j=1,2,3$ and some $\alpha_3,\alpha_6 \in (0,1)$. Let $\nu(\cdot) \in C^1(\R^M)$,
\begin{equation}
    \hat{e}_1=(1,0,...,0) \in \R^M \mbox{ \ \ and \ \ }
    \nu(\hat{e}_1) >0.
\end{equation}
Furthermore let $u_0\in \dot B^{1/2}_{2,1}\cap \dot B^{2/5}_{5/4,1}(\R^3)$,  $\vec \rho \in L_\infty(\R^3)$ and 
$\nabla \rho_{0i}^{(1-\alpha_p)} \in L_p(\R^3)$ with $p=3,6$ and $i=1,...,M$.

Then there exists a constant $\delta>0$ such that if the initial data satisfy
\begin{equation} \label{0:5}
\|u_0\|_{\dot B^{1/2}_{2,1}\cap \dot B^{2/5}_{5/4,1}}+\|\vec \rho_0-\hat{e}_1 \|_{L_\infty}
+\sum_{i=1}^M\left(\|\nabla \rho_{0i}^{(1-\alpha_3)}\|_{L_3} + \|\nabla \rho_{0i}^{(1-\alpha_3)}\|_{L_6}\right) 
\leq \delta,     
\end{equation}
then problem \eqref{sys}-\eqref{init} admits a unique global in time regular solution satisfying the estimate 
\begin{equation} \label{0:6}
\begin{aligned}
&\|u\|_{W^{2,1}_{2,(4/3,1)}\cap W^{2,1}_{5/4,(5/4,1)} }+\|tu\|_{W^{2,1}_{6,(4,1)}\cap  W^{2,1}_{2,(4,1)}
}+\|\nabla u\|_{L_1(\R_+;L_\infty)} 
\\[3pt]
&
+ \|\vec \rho - \hat{e}_1\|_{\linf} + \|\nabla \rho\|_{L_\infty(\R_+;L_3\cap L_6)} \\[3pt]
&\leq \phi( \delta ),
\end{aligned}
\end{equation}
where the space $W^{2,1}_{p,(q,r)}$ is defined in \eqref{def:3}.
\end{thm}

The motivation of the above theorem is an elegant application of the framework of the Lorentz spaces in time. The key point of the analysis of systems of the type of (\ref{sys}) is the control of 
the gradients of densities. It requires a crucial condition 
\begin{equation} \label{int:finite}
\int_0^\infty \|\nabla u(t)\|_{L_\infty} dt < +\infty.
\end{equation}
The estimates in Lorentz spaces for the momentum equation
together with the structure condition given by (\ref{omega:pos}) enable us to guarantee \eqref{int:finite}. An example of a realization of this condition is shown in Section \ref{sec:density}. This example can be generalized to a large class of production functions $\omega_i(\cdot)$. 

In the whole space \eqref{int:finite} is a trouble maker. We could obtain it using the standard $L_p-L_q$ estimates giving the bounds in time weigths $t^{-\alpha}$. However then the functional setting would be very complex, and far away from the critical and homogeneous case. The Lorentz spaces allow to work with just one time weight $t^{-1}$. Having the index 1, i.e. considering spaces $L_{p,1}$, it allows to have
$$
\int_0^\infty t^{-1/2} f(t) dt \leq C\|t^{-1/2}\|_{L_{2,\infty} }\|f\|_{L_{2,1}}.
$$
As a consequence we are allowed to work in the framework of homogeneous spaces, which makes our analysis much easier from the computational viewpoint. Nevertheless the relatively low dimension requires an extra regularity for the initial datum for the velocity
in the $\dot B^{2/5}_{5/4,1}$ space. Thanks to that we are able to control the time weighted norm in low integrability. 
Note also that the space $W^{2,1}_{5/4}$ is related to regularity of weak solutions to the Navier-Stokes equations, since for regularity coming for  weak solutions we have $u\cdot \nabla u \in L_{10/3} \cdot L_2  \subset L_{5/4}$. Therefore the regularity theory for the Stokes system implies that solutions belong to $W^{2,1}_{5/4}$. We shall underline this assumption defines the sub-critical regularity.

As a final remark we shall note that our non-standard spaces $L_{p,1}(0,T;L_q)$
can not be replaced by $L_{p,1}(\R^3 \times (0,T))$ or by $L_p(0,T;L_{q,1})$. The analysis of possible application of such setting gives a negative answer, it comes from the fact that the space $L_{q,1}(\R^3)$ is not a UMD one, in particular, it is not reflexive, so one can not apply the maximal regularity theory for the abstract semigroups, since it is based on the properties of the UMD spaces \cite{Weis,DHP}. The latter is necessary to obtain the regularity of the linearized problem in Lorentz spaces using interpolation. 




\smallskip 

The rest of the paper is organized as follows. In section \ref{sec:lin}
we show some auxiliary results and prove an estimate in Lorentz spaces for the Stokes system which is a linearization of the fluid part of system \eqref{sys}. Section \ref{sec:nonlin} is devoted to {\em a priori} estimates for system \eqref{sys}. First we prove the estimates for the velocity assuming appropriate regularity and smallness of fractional densitities. Here we follow the ideas from \cite{DMT}, hence we show the estimate for the velocity in $W^{2,1}_{2,(4/3,1)}$. Next we show a time weighted estimate, i.e. estimate for $tu$ in $W^{2,1}_{6,(4,1)}$. This allows to show that $\nabla u \in L_1(\R_+;L_\infty)$. The latter estimate allows us to show the bounds for the gradients of fractional densities, however for this purpose we need appropriate structure of functions describing the chemical reactions, which at this stage is assumed. In subsection \ref{sec:density} we derive necessary estimates of fractional densities under the assumption \eqref{omega:pos} and provide an example of realization of this assumption.
Finally in section \ref{sec:existence} we apply the estimates from section \ref{sec:nonlin} to prove the existence. The hyperbolicity of the transport equations of the species subsystem
does not allow to use the contraction procedure directly. For this reason we apply a modification of the Banach Fixed Point Theorem showing the contraction in a weaker norm than the regularity of the solution.

\section{Auxiliary results and linear theory} \label{sec:lin}

This part is dedicated to the main tool of our technique, the estimate in the Lorentz spaces adapted from \cite{DMT}. In addition we show some key estimates on products of functions.

\smallskip 

Consider the Stokes system
\begin{equation} \label{Stokes}
\begin{aligned}
& u_t-\nu \Delta  u + \nabla \pi =f, \qquad {\rm div\,} u =0 \qquad \mbox{ in } \R^3 \times [0,T),\\[5pt]
& u|_{t=0} = u_0 \qquad \mbox{ in } \R^3. 
\end{aligned}
\end{equation}
We show the following estimates for \eqref{Stokes} in Lorentz spaces: 

\begin{thm} \label{p1}
Let $T\in (0,\infty]$, $p,q \in (1,\infty)$, $r \in [1,\infty]$.
Assume $f \in L_{q,r}(0,T;L_p(\R^3))$ and $u_0 \in B^{2-2/q}_{p,r}(\R^3)$. Then \eqref{Stokes} admits a unique solution $(u,\pi)$ such that
\begin{equation} \label{est:lin1}
 \|u\|_{L_\infty(0,T;\dot B^{2-2/q}_{p,r})} + 
 \|u_t,\nu \nabla^2 u\|_{L_{q,r}(0,T;L_p)} \leq \\
 C\big( \|f\|_{L_{q,r}(0,T;L_p)} + \|u_0\|_{\dot B^{2-2/q}_{p,r}}\big),
\end{equation}
where constant $C$ in (\ref{est:lin1}) does not depend on $T$.

Moreover, if $\frac{2}{q}+\frac{d}{p}>2$ then for all $q<s<\infty$ and $p\leq m$ such that $1+\frac{d}{2}(\frac{1}{m}-\frac{1}{p})>0$ we have
\begin{equation} \label{est:lin2}
 \|u\|_{L_{s,r}(0,T;L_m)}\leq C\big( \|u\|_{L_\infty(0,T;\dot B^{2-2/q}_{p,r})} +
 \|u_t,\nabla^2u\|_{L_{q,r}(0,T;L_p)}\big)
\end{equation}
with
\begin{equation*}
 \frac{d}{2m}+\frac{1}{s}=\frac{1}{q}+\frac{d}{2p}-1.
\end{equation*}
And, if $\frac{2}{q}+\frac{d}{p}>1$ then for all $q<s<\infty$ and $p\leq m$ such that $1+{d}(\frac{1}{m}-\frac{1}{p})>0$ we have 
\begin{equation} \label{est:lin3}
 \|\nabla u\|_{L_{s,r}(0,T;L_m)}\leq C\big( \|u\|_{L_\infty(0,T;\dot B^{2-2/q}_{p,r})} +
 \|u_t,\nabla^2u\|_{L_{q,r}(0,T;L_p)}\big)
\end{equation}
with
\begin{equation*}
 \frac{d}{m}+\frac{2}{s}=\frac{2}{q}+\frac{d}{p}-1.
\end{equation*}

%
\end{thm}

{\bf Proof.} 
The existence of strong solution with estimates \eqref{est:lin1} and \eqref{est:lin2} have been shown in \cite[Proposition 2.1]{DMT} for the heat equation in the whole space. The proof is based on the maximal $L_p-L_q$ regularity for the heat equation and interpolation properties of Lorentz spaces. Proposition \ref{p1} can be proved similarly, for the sake of completeness we present the proof. We start with the following maximal regularity estimate for the Stokes problem \eqref{Stokes} in the whole space (it can be obtained for instance from Theorem 4.1 in (\cite{ShiShi}) combined with the fact that the Stokes operator generates an analytic semigroup):
\begin{equation}\label{lpstokes}
\|u\|_{L_\infty(\R_+;\dot B^{2-2/\alpha}_{p,\alpha}(\R^3))}
+\|u_t,\nabla^2 u\|_{L_\alpha(\R_+;L_p(\R^3))}
\leq C \left( \|u_0\|_{\dot B^{2-2/\alpha}_{p,\alpha}(\R^3)} + \|f\|_{L_\alpha(\R_+;L_p(\R^3))} \right) 
\end{equation}
for $1<\alpha,p<\infty$.

To start the proof we need to recall some facts from the abstract definition of the Lorentz spaces
\cite[Theorem 2:1.18.6]{Tr}. Namely, for any Banach space $A$ we have
\begin{equation}
    (L_{q_0,r_0}(0,T;A);L_{q_1,r_1}(0,T;A))_{\theta,r}=L_{q,r}(0,T;A) \mbox{ with } 
    \frac{1}{q}=\frac{\theta}{q_0}+\frac{1-\theta}{q_1}, \quad \theta \in (0,1).
\end{equation}
Taking $T=\infty$, $A=L_p(\R^d)$ ($L_p$ for short), $r_0=q_0$ and $r_1=q_1$ we have
\begin{equation}\label{a7}
    (L_{q_0}(\R_+;L_p);L_{q_1}(\R_+;L_p))_{\theta,r}=L_{q,r}(\R_+;L_p).
\end{equation}
So for any $q\in (1,\infty)$ one finds $1<q_0<q<q_1<\infty$ such that 
\begin{equation}\label{a8}
    \frac{1}{q}=\frac12 \frac{1}{q_0} + \frac12 \frac{1}{q_1}
\end{equation}
and we construct the space $L_{q,r}(\R_+;L_p)$ for any $q\in (1,\infty)$ taking $\theta =\frac12$ in (\ref{a7}).
In addition to (\ref{a7}), by \eqref{0:3a} we find 
\begin{equation*}
    (\dot B^{2-2/q_0}_{p,q_0},\dot B^{2-2/q_1}_{p,q_1})_{\frac12,r} = \dot B^{2-2/q}_{p,r}.
\end{equation*}

Now we are in a position to prove the first assertion of the Theorem. Take (\ref{lpstokes}) with $q_0$ and $q_1$ as above instead of $\alpha$, 
then applying the interpolation theorem for the real interpolation $(\cdot,\cdot)_{1/2,r}$, from \eqref{a7} and \eqref{a8} we deduce \eqref{est:lin1}.
In order to show (\ref{est:lin2}) we recall that the imbedding theorem (\ref{lpstokes}) implies that
\begin{equation}\label{lpstokes-u}
\|u\|_{L_s(\R_+;L_m(\R^3))}
\leq C \left( \|u_0\|_{\dot B^{2-2/\alpha}_{p,\alpha}(\R^3)} + \|f\|_{L_\alpha(\R_+;L_p(\R^3))} \right) 
\end{equation}
for $\frac{2}{\alpha}+\frac{d}{p}>2$, $\alpha<s<\infty$   and $p\leq m< \infty$ such that $1+\frac{d}{2}(\frac{1}{m}-\frac{1}{p})>0$ and
\begin{equation*}
 \frac{d}{2m}+\frac{1}{s}=\frac{1}{\alpha}+\frac{d}{2p}-1.
\end{equation*}
Then replacing $\alpha$ with $q_0$ and $q_1$ we get
\begin{equation}\label{lpstokes-u1}
\|u\|_{L_{s,r}(\R_+;L_m(\R^3))}
\leq C \left( \|u_0\|_{\dot B^{2-2/q}_{p,r}(\R^3)} + \|f\|_{L_{q,r}(\R_+;L_p(\R^3))} \right).
\end{equation}
It remains to show (\ref{est:lin3}). For $\nabla u$ we have the following estimate
\begin{equation}\label{lpstokes-gu}
\|\nabla u\|_{L_s(\R_+;L_m(\R^3))}
\leq C \left( \|u_0\|_{\dot B^{2-2/\alpha}_{p,\alpha}(\R^3)} + \|f\|_{L_\alpha(\R_+;L_p(\R^3))} \right) 
\end{equation}
provided $\frac{2}{\alpha}+\frac{d}{p}>1$, $\alpha<s<\infty$ and $p\leq m < \infty$ such that $1+{d}(\frac{1}{m}-\frac{1}{p})>0$ and
\begin{equation*}
 \frac{d}{m}+\frac{2}{s}=\frac{2}{\alpha}+\frac{d}{p}-1.
\end{equation*}
Then again using the interpolation theory we get
\begin{equation}\label{lpstokes-gu}
\|\nabla u\|_{L_{s,r}(\R_+;L_m(\R^3))}
\leq C \left( \|u_0\|_{\dot B^{2-2/q}_{p,r}(\R^3)} + \|f\|_{L_{q,r}(\R_+;L_p(\R^3))} \right),
\end{equation}
as
\begin{equation*}
 \frac{d}{m}+\frac{2}{s}=\frac{2}{q}+\frac{d}{p}-1.
\end{equation*}
\qed



\noindent
Let us finish this section with deducing from \eqref{holder} the following useful inequalities:

\begin{lem} For any $\Omega \subset \R^n$, $n \in \N$ and $1<q<\infty$ we have
\begin{align}
& \|fg\|_{L_{q,1}(\R_+;L_2(\Omega))} \leq C \|f\|_{L_\infty(\R_+,L_3(\Omega))}\|g\|_{L_{q,1}(\R_+;L_6(\Omega))} \label{2:1}, \\
& \|fg\|_{L_{q,1}(\R_+;L_{5/4}(\Omega))} \leq C \|f\|_{L_\infty(\R_+,L_3(\Omega))}\|g\|_{L_{q,1}(\R_+;L_{15/7}(\Omega))}. \label{2:2}
\end{align}
\end{lem}

{\bf Proof.} By standard H\"older inequality and its counterpart for Lorentz spaces \eqref{holder} we have 
$$
\Big\| \|fg(\cdot)\|_{L_2(\Omega)} \Big\|_{L_{q,1}(\R_+)} \leq
\Big\| \|f(\cdot)\|_{L_3(\Omega)}\|g(\cdot)\|_{L_6(\Omega)} \Big\|_{L_{q,1}(\R_+)} \leq \|f\|_{L_\infty(\R_+;L_3(\Omega))}\|g\|_{L_{q,1}(\R_+;L_6(\Omega))}
$$
which proves \eqref{2:1}. \eqref{2:2} is shown similarly using 
$\|fg\|_{L_{5/4}(\Omega)}\leq \|f\|_{L_3(\Omega)}\|g\|_{L_{15/7}(\Omega)}$.

\begin{flushright}
$\square$
\end{flushright}

\section{A priori estimates} \label{sec:nonlin}

We start our investigation of system (\ref{sys}) with finding the a priori bounds of solutions (\ref{0:6}),  provided the data are sufficiently small. It will be split into several steps stated as lemmas.



\subsection{Velocity bounds}

In this section we show the required a priori estimates for the velocity. At this moment we assume that the whole $\vec \rho$ is known and it is sufficiently close to constant vector. 
The structure of the system is similar to the classical Navier-Stokes system, the coupling with reaction equations is via the 
viscosity coefficient. In the framework of small solutions one can write the momentum equation in the following way
\begin{equation} \label{stokes1}
 \begin{array}{lcl}
 u_t  - \bar \nu \Delta u + \nabla p = F  & \mbox{ in } & \R^3 \times [0,T),\\[5pt]
 \div u=0 & \mbox{ in } & \R^3 \times [0,T),\\[5pt]
 u|_{t=0}=u_0,
 \end{array}
\end{equation}
with 
\begin{equation} \label{def:F}
F=(1-\rho) u_t   -  \rho u\cdot \nabla u +(\nu(\vec \rho\,) -\bar \nu) \Delta u + 
 \nabla_{\vec \rho} \, \nu(\vec \rho\,) \nabla \rho : \D(u), 
\end{equation}
where 
\begin{equation}\label{e1}
  \quad \bar\nu = \nu(\hat{e}_1) \mbox{ \ \ with \ \ }  \hat{e}_1=(1,0,...,0) \in \R^M.
\end{equation}

The first step is to show the following bound. The below regularity, provided the density is a given sufficiently smooth function, guarantees the unique solvability of the Navier-Stokes equations \cite{Sohr}, by the imbedding $\dot B^{1/2}_{2,1}(\R^3) \subset L_3(\R^3)$.
%
\begin{lem} \label{l:2}
%
Assume $u$ solves \eqref{stokes1}-\eqref{def:F}. Then we have the following inequality 
\begin{equation} \label{1:10}
\begin{aligned}
& \|u\|_{W^{2,1}_{2,(4/3,1)}} :=\|u\|_{L_\infty(\R_+;\dot B^{1/2}_{2,1})} + \|u_t,\nabla^2 u \|_{L_{4/3,1}(\R_+;L_2)}\\[3pt]
&\lesssim (\|\vec \rho - \hat{e}_1\|_{L_\infty(\R_+;L_\infty)}+
 \|\nabla \vec \rho \,\|_{L_\infty(\R_+;L_3)}+
 \|u\|_{W^{2,1}_{2,(4/3,1)}})\|u\|_{W^{2,1}_{2,(4/3,1)}}
 + \|u_0\|_{\dot B^{1/2}_{2,1}}.
 \end{aligned}
\end{equation}
\end{lem}
{\bf Proof.} We shall apply Theorem \ref{p1}, therefore we have to estimate step by step the terms from $F$ in the required regularity. First, 
$F_1=(1-\rho) u_t $ gives
\begin{equation*}
    \|F_1\|_{L_{4/3,1}(\R_+;L_2)} \leq C\|1-\rho\|_{L_\infty(\R_+;L_\infty)}
    \|u_t\|_{L_{4/3,1}(\R_+;L_2)}.
\end{equation*}
Next, for $F_2=\rho u\cdot \nabla u$ we have by \eqref{2:1} with $q=\frac{4}{3}$ 
\begin{equation*}
\|F_2\|_{L_{4/3,1}(\R_+;L_2)} \leq C \|u \cdot \nabla u \|_{L_{4/3,1}(\R_+;L_2)}\leq C\|u\|_{L_\infty(\R_+;L_3)} 
 \|\nabla u \|_{L_{4/3,1}(\R_+;L_6)}.
\end{equation*}
For $F_3= (\nu(\vec \rho\,) -\bar \nu) \Delta u$ we first notice that 
\begin{equation} \label{visc}
\|\nu(\vec\rho\,) - \bar \nu\|_{L_\infty(\R_+;L_\infty)}\leq \|\vec\rho\ -\hat{e}_1\|_{L_\infty(\R_+;L_\infty)},
\end{equation}
therefore
\begin{equation*}
    \|F_3\|_{L_{4/3,1}(\R_+;L_2)} \leq 
    C\|\vec \rho -\vec \rho_0\|_{L_\infty(\R_+;L_\infty)}
    \|\nabla^2 u\|_{L_{4/3,1}(\R_+;L_2)}.
\end{equation*}
The last term $F_4= \nabla_{\vec \rho}\, \nu(\vec \rho\,) \nabla \vec \rho : \D(u)$ coming from the variability of $\nu(\cdot)$ can be bounded again using \eqref{2:1}:
\begin{equation*}
\|F_4\|_{L_{4/3,1}(\R_+;L_2)}=\| \nabla_{\vec \rho} \, \nu(\vec \rho\,) \cdot \nabla \vec \rho : \D(u)\|_{L_{4/3,1}(\R_+;L_2)}\leq C\|\nabla \vec \rho\, \|_{L_\infty(\R_+;L_3)} \|\nabla u\|_{L_{4/3,1}(\R_+;L_6)}.
\end{equation*}
Combining the estimates for $F_1-F_4$ with Theorem \ref{p1} we obtain 
\begin{equation*} \label{1:15}
\begin{aligned}
&   \|u\|_{L_\infty(\R_+;\dot B^{1/2}_{2,1})} + \|u_t,\nabla^2 u \|_{L_{4/3,1}(\R_+;L_2)}\leq 
 C( \|F_1,F_2,F_3,F_4\|_{L_{4/3,1}(\R_+;L_2)} + \|u_0\|_{\dot B^{1/2}_{2,1}})\\
& \leq C\Big[ (   \|\vec \rho - \hat{e}_1\|_{L_\infty(\R_+;L_\infty)} 
  + \|\nabla \vec \rho\,\|_{L_\infty(\R_+;L_3)}
  +\|u\|_{L_\infty(\R_+;L_3)})  
\|u_t,\nabla^2 u \|_{L_{4/3,1}(\R_+;L_2)} 
+ \|u_0\|_{\dot B^{1/2}_{2,1}}\Big].
\end{aligned}
\end{equation*}
Since 
\begin{equation} \label{imbed:3}
\dot B^{1/2}_{2,1}(\R^3) \subset L_3(\R^3), 
\end{equation}
we obtain the desired bound (\ref{1:10}).
\qed

\begin{cor}
Note that the imbeddings \eqref{est:lin2}-\eqref{est:lin3} imply that
\begin{equation}\label{paraimb}
    \|u\|_{L_{4,1}(\R_+;L_6)} +
    \|\nabla u\|_{L_{4,1}(\R_+;L_2)}\leq  C \|u\|_{W^{2,1}_{2,(4/3,1)}}.
\end{equation}
\end{cor}

The second step is to show the  bound in the lower regularity, it appears that  bounds coming from (\ref{paraimb}) are too high for dimension three.
%

\begin{lem} \label{l:2a}
Assume $u$ solves \eqref{stokes1}-\eqref{def:F}
with $u_0 \in \dot B^{2/5}_{5/4,1}(\R^3)$.
Then the following inequality  is valid
\begin{equation} \label{1:10a}
\begin{aligned}
& \|u\|_{W^{2,1}_{5/4,(5/4,1)}} :=\|u\|_{L_\infty(\R_+;\dot B^{2/5}_{5/4,1})} + \|u_t,\nabla^2 u \|_{L_{5/4,1}(\R_+;L_{5/4})}\\[3pt]
 &\lesssim (\|\vec \rho - \hat{e}_1\|_{L_\infty(\R_+;L_\infty)}+
 \|\nabla \vec \rho \,\|_{L_\infty(\R_+;L_3)}+
 \|u\|_{W^{2,1}_{2,(4/3,1)}})\|u\|_{W^{2,1}_{5/4,(5/4,1)}}
 + \|u_0\|_{\dot B^{2/5}_{5/4,1}}.
 \end{aligned}
\end{equation}
\end{lem}
{\bf Proof.} We proceed similarly as in the proof of Lemma \ref{l:2}. Firstly,
\begin{equation*}
    \|F_1\|_{L_{5/4,1}(\R_+;L_{5/4})} \leq C\|1-\rho\|_{L_\infty(\R_+;L_\infty)}
    \|u_t\|_{L_{5/4,1}(\R_+;L_{5/4})}.
\end{equation*}
Next, for $F_2$ we have by \eqref{2:2} with $q=\frac{5}{4}$ 
\begin{equation*}
\|F_2\|_{L_{5/4,1}(\R_+;L_{5/4})} \leq C \|u \cdot \nabla u \|_{L_{5/4,1}(\R_+;L_{5/4})}\leq C\|u\|_{L_\infty(\R_+;L_3)} 
 \|\nabla u \|_{L_{5/4,1}(\R_+;L_{15/7})}.
\end{equation*}
For $F_3$, by \eqref{visc} 
\begin{equation*}
    \|F_3\|_{L_{5/4,1}(\R_+;L_{5/4})} \leq 
    C\|\vec \rho -\vec \rho_0\|_{L_\infty(\R_+;L_\infty)}
    \|\nabla^2 u\|_{L_{5/4,1}(\R_+;L_{5/4})}.
\end{equation*}
The last term $F_4$ can be bounded again using \eqref{2:2}
with $q=\frac{5}{4}$
\begin{equation*}
\|F_4\|_{L_{5/4,1}(\R_+;L_{5/4})}=\| \nabla_{\vec \rho} \, \nu(\vec \rho\,) \cdot \nabla \vec \rho : \D(u)\|_{L_{5/4,1}(\R_+;L_{5/4})}\leq C\|\nabla \vec \rho\, \|_{L_\infty(\R_+;L_3)} \|\nabla u\|_{L_{5/4,1}(\R_+;L_{15/7})}.
\end{equation*}
Combining the estimates for $F_1-F_4$ with Theorem \ref{p1} and the imbedding ($\dot W^1_{5/4}(\R^3) \subset L_{15/7}(\R^3)$) we obtain 
\begin{equation*} \label{1:15}
\begin{aligned}
&   \|u\|_{L_\infty(\R_+;\dot B^{2/5}_{5/4,1})} + \|u_t,\nabla^2 u \|_{L_{5/4,1}(\R_+;L_{5/4})}\leq
 C( \|F_1,F_2,F_3,F_4\|_{L_{5/4,1}(\R_+;L_{5/4})} + \|u_0\|_{\dot B^{2/5}_{5/4,1}})\\
& \leq C\Big[ (   \|\vec \rho - \hat{e}_1\|_{L_\infty(\R_+;L_\infty)} 
  + \|\nabla \vec \rho\,\|_{L_\infty(\R_+;L_3)}
  +\|u\|_{L_\infty(\R_+;L_3)})  
\|u_t,\nabla^2 u \|_{L_{5/4,1}(\R_+;L_{5/4})} 
+ \|u_0\|_{\dot B^{2/5}_{5/4,1}}\Big],
\end{aligned}
\end{equation*}
and by \eqref{imbed:3} we get \eqref{1:10a}.
\qed
\begin{cor}
Note that the imbeddings \eqref{est:lin2} imply that
\begin{equation}\label{paraimb-a}
    \|u\|_{L_{4,1}(\R_+;L_2)} 
    \leq  C \|u\|_{W^{2,1}_{5/4,(5/4,1)}}.
\end{equation}
\end{cor}


The next step is to improve the regularity without changing the initial data. We  consider  time weighted spaces, characteristic for the approach via the semigroup theory. 
Hence we aim at showing the following

\begin{lem} \label{l:21}
Assume $u$ solves \eqref{stokes1}-\eqref{def:F}. Then
\begin{equation} \label{1:17}
\begin{aligned}
 \|tu\|_{W^{2,1}_{2,(4,1)}} \lesssim &
 \Big(\|\vec \rho-\hat{e}_1\|_{L_\infty(\R^3 \times \R_+)} +\|\nabla \vec\rho\, \|_{L_\infty(\R_+;L_3)} + 
 \|u\|_{W^{2,1}_{2,(4/3,1)}} \Big)\|tu\|_{W^{2,1}_{2,(4,1)}}\\
&+ \|u\|_{W^{2,1}_{5/4,(5/4,1)}}.
 %
\end{aligned}
\end{equation}
\end{lem}
{\bf Proof.}
We multiply \eqref{stokes1}$_1$-\eqref{def:F} by $t$ obtaining 
\begin{equation} \label{subsys:t}
 \begin{array}{lcl}
 (tu)_t  - \bar \nu \Delta (tu) + \nabla (tp) = F^t  & \mbox{ in } & \R^3 \times [0,T),\\[5pt]
 \div (tu)=0 & \mbox{ in } & \R^3 \times [0,T)
 \end{array}
\end{equation}
with 
\begin{equation} \label{def:Ft}
F^t=(1-\rho) (tu)_t + \rho u  -  t\rho u\cdot \nabla u +(\nu(\vec \rho) -\bar \nu) \Delta (tu) + {
 t\nabla_{\vec \rho} \, \nu(\vec \rho\,) \nabla \rho : \D(u)}.
\end{equation}
In order to show \eqref{1:17} we apply Theorem \ref{p1} to \eqref{subsys:t}. Therefore we estimate $\|F^t\|_{L_{4,1}(\R_+;L_2)}$. 
We have 
\begin{equation} \label{1:18a}
\|(1-\rho)(tu)_t\|_{L_{4,1}(\R_+;L_2)} \leq C \|(1-\rho)\|_{L_\infty(\R_+;L_\infty)} \|(tu)_t\|_{L_{4,1}(\R_+;L_2)} \end{equation}
and, since we assume the smallness of perturbation of the density, by \eqref{paraimb} we obtain
\begin{equation} \label{1:18b}
\|\rho u\| _{L_{4,1}(\R_+;L_2)} \lesssim \|u\|_{L_{4,1}(\R_+;L_2)} 
\lesssim  \|u\|_{W^{2,1}_{5/4,(5/4,1)}}.
\end{equation}
Also, similarly as in the proof of \eqref{1:10},
\begin{equation} \label{1:18c}
    \|(\nu(\vec \rho\,) -\bar \nu) \Delta (tu)\|_{L_{4,1}(\R_+;L_2)} \lesssim 
    \|\vec \rho -\hat{e}_1\|_{L_\infty(\R_+;L_\infty)}
    \|\nabla^2 (tu)\|_{L_{4,1}(\R_+;L_2)}.
\end{equation}
Next we  verify if 
\begin{equation}\label{1:19}
 tu \nabla u \in L_{4,1}(\R_+;L_2).
\end{equation}
For this purpose, by \eqref{2:1} with $q=4$ we observe that
\begin{equation} \label{1:23}
    \|tu \nabla u \|_{L_{4,1}(\R_+;L_2)}\leq 
    \|u\|_{L_\infty(\R_+;L_3)}\|t\nabla u \|_{L_{4,1}(\R_+;L_6)}.
\end{equation}
In order to complete the bound for $F^t$ we estimate the last term. For this purpose we apply again \eqref{2:1} with $q=4$ to get 
\begin{equation} \label{1:24}
\begin{aligned}
 \| \nabla_{\vec \rho}\, \nu(\vec \rho\,) \nabla \vec \rho :t \D(u)
 \|_{L_{4,1}(\R_+;L_2)} & \leq
 C\|\nabla \vec \rho \,\|_{L_\infty(\R_+;L_3)} \| t \nabla u\|_{L_{4,1}(\R_+;L_6)}\\[3pt]
& \leq C\|\nabla \vec \rho\,\|_{L_\infty(\R_+;L_3)} \|tu\|_{W^{2,1}_{2,(4,1)}}.
\end{aligned}
\end{equation}
Combining \eqref{1:18a}-\eqref{1:18c}, \eqref{1:23} and \eqref{1:24} we get 
\begin{equation*} \label{1:25}
\begin{aligned}
\|F^t\|_{L_{4,1}(\R_+;L_2)} \lesssim & \Big(\|\vec \rho-\hat{e}_1\|_{L_\infty(\R^3 \times \R_+)} +\|\nabla \vec\rho\, \|_{L_\infty(\R_+;L_3)} + \|u\|_{L_\infty(\R_+;L_3)} \Big)\|tu\|_{W^{2,1}_{2,(4,1)}}
\\ & + \|u\|_{W^{2,1}_{5/4,(5/4,1)}}.
\end{aligned}
\end{equation*}
%
\qed

\begin{lem} \label{l:21:2}
Assume $u$ solves \eqref{stokes1}-\eqref{def:F}. Then
\begin{equation} \label{1:17:2}
\begin{aligned} 
 \|tu\|_{W^{2,1}_{6,(4,1)}} \lesssim 
& \Big(\|\vec \rho-\hat{e}_1\|_{L_\infty(\R^3 \times \R_+)} +\|\nabla \vec\rho\, \|_{L_\infty(\R_+;L_6)} + 
 \|u\|_{W^{2,1}_{2,(4/3,1)}} \Big)\|tu\|_{W^{2,1}_{6,(4,1)}}\\
& 
+\|u\|_{W^{2,1}_{2,(4/3,1)}} + \|\nabla \vec\rho\,\|_{L_\infty(\R_+;L_6)}\|tu\|_{W^{2,1}_{2,(4,1)}}.
 %
\end{aligned}
\end{equation}
\end{lem}

{\bf Proof.}
%
We proceed similarly as in the proof of Lemma \ref{l:21}, applying Theorem \ref{p1} to \eqref{subsys:t}. For this purpose we estimate $\|F^t\|_{L_{4,1}(\R_+;L_6)}$. 
We have 
\begin{equation} \label{1:18a:2}
\|(1-\rho)(tu)_t\|_{L_{4,1}(\R_+;L_6)} \leq C \|(1-\rho)\|_{L_\infty(\R_+;L_\infty)} \|(tu)_t\|_{L_{4,1}(\R_+;L_6)} \end{equation}
and 
\begin{equation} \label{1:18b:2}
\|\rho u\| _{L_{4,1}(\R_+;L_6)} \lesssim \|u\|_{L_{4,1}(\R_+;L_6)} 
\lesssim  \|u\|_{W^{2,1}_{2,(4/3,1)}}.
\end{equation}
Also, similarly as before, 
\begin{equation} \label{1:18c:2}
\begin{aligned}
    \|(\nu(\vec \rho\,) -\bar \nu) \Delta (tu)\|_{L_{4,1}(\R_+;L_6)} 
    &\lesssim \|\vec \rho -\hat{e}_1\|_{L_\infty(\R_+;L_\infty)}
    \|tu\|_{W^{2,1}_{6,(4,1)}}.
\end{aligned}
\end{equation}
Two remaining terms in $F^t$ are more demanding.
First we need to verify if 
\begin{equation}\label{1:19:2}
 tu \nabla u \in L_{4,1}(\R_+;L_6).
\end{equation}
Let us start with observing that 
\begin{equation} \label{1:20:2}
L_2(\R^3) \subset \dot W^{-1}_{6}(\R^3).    
\end{equation}
In order to show \eqref{1:20:2} recall that $\dot W^{-1}_6(\R^3)$ is a space of linear functionals over the homogeneous space $\dot W^1_{6/5}(\R^3)$, so $f\in \dot W^{-1}_6(\R^3)$ iff for all $g \in \dot W^1_{6/5}(\R^3)$ 
there exists a constant $C$ such that
\begin{equation*}
 \left| \int_{\R^3} fg dx\right| \leq
 C \|g\|_{\dot W^1_{6/5}}.
\end{equation*}
From the Sobolev imebeddings 
\begin{equation*}
 \dot W^1_{6/5}(\R^3) \subset L_2(\R^3).
\end{equation*}
So we  obtain \eqref{1:20:2}, since for all $f\in L_2(\R^3)$
\begin{equation*}
 \left| \int_{\R^3} fg dx\right| \leq \|f\|_{L_2} \|g\|_{L_2}\leq 
 C \|f\|_{L_2} \|g\|_{\dot W^1_{6/5}}, 
\end{equation*}
which implies $\|f\|_{\dot W^{-1}_6} \leq C\|f\|_{L_2}$. Now we are ready to show (\ref{1:19:2}). Note that by the interpolation inequality
\begin{equation*}
    \|t^{1/2} \nabla u\|_{L_6} \leq C\|t \nabla u\|_{\dot W^1_6}^{1/2}\|\nabla u\|_{\dot W^{-1}_6}^{1/2},
\end{equation*}
so by \eqref{holder} and \eqref{0:2} 
\begin{equation*}
\begin{aligned}
    \left\| \|t^{1/2} \nabla u\|_{L_6} \right\|_{L_{4,1}(\R_+)} &\leq C
    \left\| \|t \nabla u\|_{\dot W^1_6}^{1/2} \right\|_{L_{8,2}(\R_+)} \left\| \|\nabla u\|_{\dot W^{-1}_6}^{1/2}\right\|_{L_{8,2}(\R_+)}\\[3pt]
    &\leq C\|t \nabla u\|_{L_{4,1}(\R_+;\dot W^{1}_6)}^{1/2}
    \|\nabla u\|_{L_{4,1}(\R_+;\dot W^{-1}_6)}^{1/2}.
    \end{aligned}
\end{equation*}
Therefore by \eqref{1:20:2} 
\begin{equation} \label{1:21a}
\|t^{1/2}\nabla u\|_{L_{4,1}(\R_+;L_6)} \leq C\|t \nabla u\|_{L_{4,1}(\R_+;\dot W^{1}_6)}^{1/2}
    \|\nabla u\|_{L_{4,1}(\R_+;L_2)}^{1/2}.
\end{equation}
Moreover, due to \eqref{0:3} and \eqref{0:4} we find that
\begin{equation*} 
\|t^{1/2} u\|_{L_\infty} \leq C\|t^{1/2}u\|_{\dot B^{1/2}_{6,1}} \leq C\|tu\|_{\dot B^{3/2}_{6,1}}^{1/2}\|u\|_{\dot B^{-1/2}_{6,1}}^{1/2},
\end{equation*}
so
\begin{equation} \label{1:22a}
    \|t^{1/2} u\|_{L_\infty(\R_+;L_\infty)} \leq 
    C \|tu\|_{L_\infty(\R_+;\dot B^{3/2}_{6,1})}^{1/2}\|u\|_{L_\infty(\R_+;\dot B^{-1/2}_{6,1})}^{1/2}.
\end{equation}
Combining \eqref{1:21a} and \eqref{1:22a} with Young inequality we obtain 
\begin{equation} \label{1:23:2}
\begin{aligned}
\|tu\nabla u\|_{L_{4,1}(\R_+;L_6)} &\leq \|t^{1/2}u\|_{L_\infty(\R^3\times \R_+)}\|t^{1/2}\nabla u\|_{L_{4,1}(\R_+;L_6)} \\[3pt]
&\lesssim \|u\|_{L_\infty(\R_+;B^{-1/2}_{6,1})}^{1/2}\|\nabla u\|_{L_{4,1}(\R_+;L_2)}^{1/2}
\Big( \|tu\|_{L_\infty(\R_+;\dot B^{3/2}_{6,1})} + \|t\nabla u\|_{L_{4,1}(\R_+;\dot W^1_6)} \Big) \\[3pt]
&= \|u\|_{L_\infty(\R_+;B^{-1/2}_{6,1})}^{1/2}\|\nabla u\|_{L_{4,1}(\R_+;L_2)}^{1/2} \|tu\|_{W^{2,1}_{6,(4,1)}}.
\end{aligned}
\end{equation}
It remains to estimate the last term of $F^t$. Here is the part which requires the estimate from Lemma \ref{l:21}. By the interpolation inequality we have
\begin{equation*}
    \|\nabla u \|_{L_\infty(\R^3)} \leq C
    \|\nabla u\|_{L_6(\R^3)}^{1/2}
    \|\nabla^2 u\|_{L_6(\R^3)}^{1/2}.
\end{equation*}
Thus, we conclude
\begin{equation} \label{1:24:2}
\begin{aligned} 
 \| \nabla_{\vec \rho}\, \nu(\vec \rho\,) \nabla \vec \rho :t \D(u)
 \|_{L_{4,1}(\R_+;L_6)} & \leq 
 C\|\nabla \vec \rho \,\|_{L_\infty(\R_+;L_6)} \| t \nabla u\|_{L_{4,1}(\R_+;L_\infty)}\\[3pt]
& \leq C\|\nabla \vec \rho\,\|_{L_\infty(\R_+;L_6)}\|tu\|_{W^{2,1}_{2,(4,1)}}^{1/2} \|tu\|_{W^{2,1}_{6,(4,1)}}^{1/2}.
\end{aligned}
\end{equation}
Combining \eqref{1:18a:2}-\eqref{1:18c:2}, \eqref{1:23:2} and \eqref{1:24:2} we get 
\begin{equation*} \label{1:25}
\begin{aligned}
\|F^t\|_{L_{4,1}(\R_+;L_6)} \lesssim & \Big(\|\vec \rho-\hat{e}_1\|_{L_\infty(\R^3 \times \R_+)} +
\|u\|_{L_\infty(\R_+;\dot B^{-1/2}_{6,1})}^{1/2}\|\nabla u\|_{L_{4,1}(\R_+;L_2)}^{1/2} \Big)\|tu\|_{W^{2,1}_{6,(4,1)}}
 + \\ &
 \|\nabla \vec \rho\,\|_{L_\infty(\R_+;L_6)}\|tu\|_{W^{2,1}_{2,(4,1)}}^{1/2} \|tu\|_{W^{2,1}_{6,(4,1)}}^{1/2}
 + \|u\|_{W^{2,1}_{2,(4/3,1)}}.
\end{aligned}
\end{equation*}
By (\ref{paraimb}) and imbedding $\dot B^{1/2}_{2,1}(\R^3) \subset \dot B^{-1/2}_{6,1}(\R^3)$, applying Young inequality to the first term of the second line we conclude (\ref{1:17}).
\qed

To finish the considerations for the velocity we prove the following imebedding. Let us note that this result is the heart of out technique based on the Lorentz spaces since it allows to define globally the characteristics for the species subsystem.
\begin{lem}\label{lem:4}
Let $u\in W^{2,1}_{2,(4/3,1)}(\R^3 \times \R_+)$ and 
$tu \in W^{2,1}_{6,(4,1)}(\R^3 \times \R_+)$, then $\nabla u \in L_1(\R_+;L_\infty)$ and
\begin{equation*}
    \int_0^\infty \|\nabla u(t)\|_{L_\infty} dt \leq C
    \|tu \|_{W^{2,1}_{6,(4,1)}}^{1/2} \|u\|_{W^{2,1}_{2,(4/3,1)}}^{1/2}.
    \end{equation*}
\end{lem}

{\bf Proof.} 
Note that by \eqref{0:3} and Sobolev imbedding we have
\begin{equation*}
 \|\nabla u \|_{\infty} \lesssim  \|\nabla u\|_{\dot B^{1/2}_{6,1}} \lesssim 
 \|\nabla u\|^{1/2}_{\dot W^1_6} \|\nabla u\|^{1/2}_{L_6} \lesssim \|\nabla u\|^{1/2}_{\dot W^1_6} \|\nabla u\|^{1/2}_{\dot W^1_2},
\end{equation*}
therefore by \eqref{holder} 
\begin{equation*} \label{1:30}
\begin{aligned}
 \int_0^\infty \|\nabla u\|_{\infty} dt & \leq \int_0^\infty t^{-1/2} \|t \nabla u\|^{1/2}_{W^1_6} \|\nabla u \|_{W^1_2}^{1/2}dt = \Big\|\, t^{-1/2} \|t \nabla u\|^{1/2}_{W^1_6} \|\nabla u \|_{W^1_2}^{1/2}\, \Big\|_{L_{1,1}(\R_+)}  \\
& \leq \|t^{-1/2}\|_{L_{2,\infty}(\R_+)} \left\| \|t \nabla u\|^{1/2}_{W^1_6} \right\|_{L_{8,2}(\R_+)}
 \left\| \|\nabla u\|^{1/2}_{W^1_2} \right\|_{L_{8/3,2}(\R_+)} \\
& \leq \|t^{-1/2}\|_{L_{2,\infty}(\R_+)} \left\| \|t \nabla u\|_{W^1_6} \right\|_{L_{4,1}(\R_+)}^{1/2}
 \left\| \|\nabla u\|_{W^1_2} \right\|_{L_{4/3,1}(\R_+)}^{1/2}
 \\
& \leq C\|tu\|_{L_{4,1}(0,\infty;\dot W^2_6)}^{1/2} \|u\|_{L_{4/3,1}(0,\infty;\dot W^2_{2})}^{1/2}.
 \end{aligned}
\end{equation*}

\qed

\subsection{Estimates for the density} \label{sec:density}

In this section we show the $L_p$ estimates for $\eqref{subsys}$ with given velocity field $u$ provided appropriate structural assumptions on the chemical production rates. In view of the velocity estimates we need to estimate gradients of fractional densities in $L_3 \cap L_6$, but for the sake of completness we show the estimates for $p>2$ under more general assumptions than \eqref{omega:pos}.

The following lemma delivers the basic point-wise bounds for species' densities. 
\begin{lem} \label{l:3}
Let $a_1,\ldots,a_k$, $b_1,\ldots,b_l$ and $w$ be sufficiently smooth solutions to system (\ref{subsys}) and their initial data $a^0_1,\ldots,a^0_k$, $b^0_1,\ldots,b^0_l$ and $w^0$ be nonnegative, then
they remain nonnegative for all times and the following estimates hold
\begin{equation}\label{est-rho}
    \begin{aligned}
&\displaystyle 
\sup_{t\in \R_+} \|w(t)\|_{L_\infty} \leq \|w_0\|_{L_\infty}, \qquad
\sup_{t\in \R_+} \|a_i(t)\|_{L_\infty} \leq \|a^0_j\|_{L_\infty}, \quad i=1,\ldots, k, \\[4pt]
&\displaystyle \sup_{t\in \R_+} \|b_j(t)\|_{L_\infty} \leq 
\|b_j^0\|_{L_\infty}+\theta_j \|\sum_{i=1}^k a_i^0\|_{L_\infty}, \quad j=1,\ldots,l.
    \end{aligned}
\end{equation}
\end{lem}

{\bf Proof.} The bound for $w$ is classical,
it obeys the pure transport equation. For $a_i$ note that 
by \eqref{omega:zero} $a_i(t)$ stays nonnegative as the initial datum is nonnegative and then
the maximum principle implies the desired bound.

The case of $b_j(t)$ is more involved. The non-negativity is clear since the RHS is non-negative. But we note that
$$
\partial_t( \sum_{j=1}^l b_j + \sum_{i=1}^k
a_i)
+u\cdot \nabla ( \sum_{j=1}^l b_j + \sum_{i=1}^ka_i)=0.
$$
It leads to a simple but rough bound
\begin{equation*}
    \sup_t \|b_j(t)\|_{L_\infty}
    \leq \sup_t \|\sum_{j=1}^l b_j(t) + \sum_{i=1}^ka_i(t)\|_{L_\infty} \leq
    \|\sum_{j=1}^l b_j^0 + \sum_{i=1}^ka_i^0\|_{L_\infty}.
\end{equation*}
But we can improve it. Looking at the equation for $b_j$, i.e. $(\ref{subsys})_2$,
we split 
\begin{equation}\label{b-split}
b_j(t)=b_j^{ini}(t) + b_j^{zero}(t),
\end{equation}
where 
\begin{equation} \label{b-ini}
    \partial_tb_j^{ini} + u\cdot \nabla 
    b_j^{ini}=0  \mbox{ \ \ with \ \ }
    b_j^{ini}|_{t=0}=b_j^0
\end{equation}
and
\begin{equation}\label{b-zero}
    \partial_t b_j^{zero} + u \cdot \nabla b_j^{zero}= \theta_j \sum_{i=1}^k \omega_i(\vec a) \mbox{ \ \ with \ \ }
    b_j^{zero}|_{t=0}=0.
\end{equation}
The form of (\ref{b-zero}) allows to introduce a non-negative function 
$B=\sum_{j=1}^lb_j^{zero}$ which satisfies
(recall that $\sum_{j=1}^l\theta_j=1$)
\begin{equation}\label{B-zero}
    B_t + u\cdot \nabla B = \sum_{i=1}^k \omega_i(\vec a)
    \mbox{ \ \ with \ \ } B|_{t=0} =0.
\end{equation}
Combining \eqref{B-zero} with all equations for $\vec a$ we easily deduce
\begin{equation}\label{est-B}
    \sup_t \|B(t)\|_{L_\infty}\leq \sup_t \|B(t) +
    \sum_{i=1}^k a_i(t)\|_{L_\infty}\leq \|\sum_{i=1}^k a_i^0\|_{L_\infty}.
\end{equation}
Next we observe that the regularity of $u$ implies the uniqueness for 
(\ref{b-zero}) and (\ref{B-zero}), therefore we find the following relation
\begin{equation}\label{b-0-B}
    b_j^{zero}(t)=\theta_j B(t).
\end{equation}
Taking into account the above trick we find that
\begin{equation*}
    \sup_t\|b_j(t)\|_{L_\infty} \leq 
    \|b_j^0\|_{L_\infty}  + \sup_t \theta_j\|B(t)\|_{L_\infty},
    \end{equation*}
and by (\ref{est-B}) we are done.
\qed 
The next issue is to find a good bound on the gradient of $\vec \rho$. 
In view of the assumptions of Lemmas \ref{l:2} and \ref{l:3} we need to estimate $\|\nabla \vec\rho\,\|_{L_p}$ for $p=3,6$  in terms of the velocity, but, as already mentioned, we show a general estimate in $L_p$. 

\begin{lem} \label{l:4}
Assume $\nabla u \in L_1(0,T;L_\infty)$ for some $T \in (0,+\infty]$ and $\nabla (\vec \rho\,^0_i)^{(1-\alpha_p)} \in L_p(\R^3)$ for some  $2 < p <\infty$. Assume moreover that
\begin{equation}
%
\sum_{m=1}^k a_m^{-(p-1)\alpha_p}|a_{m,x_j}|^{p-2}a_{m,x_j} \de_{x_j} \big( \omega_m(\vec a)a_m^{-\alpha_p} \big) \geq 0.
 \label{omega:pos2}    
\end{equation}
Then the solution to \eqref{subsys} satisfy 
\begin{equation} \label{3:6}
\sup_{t\in(0,T)} \|\nabla \vec\rho(t)\|_{L_p} \leq 
\phi\Big(\|\nabla \vec \rho\,^0\|_{L_p},\sum_{i=1}^k\|\nabla (a_i^0)^{1-\alpha_p}\|_{L_p}\Big)\exp\left(C \int_0^\infty \|\nabla u(s)\|_{L_\infty}ds \right),   
\end{equation}
where $\phi$ is a continuous, positive function such that $\phi(0)=0$ and $C=C(p,\alpha_p)$.
\end{lem}

{\bf Proof.} Note first that as $\nabla \rho_i^{(1-\alpha)} \in L_p(\R^3)$ and 
$\rho_i^0$ is small, then $\nabla \rho^0_i\in L_p(\R^3)$, and it is small too. Indeed, we have
\begin{equation}
    \nabla \rho_i = \frac{1}{1-\alpha} \rho_i^{\alpha} \nabla \rho_{i}^{(1-\alpha)},
\end{equation}
thus
\begin{equation}\label{yy1}
    \|\nabla \rho_i\|_{L_p}\leq C_\alpha \|\rho_i\|_{L_\infty}^\alpha \|\nabla \rho_i^{(1-\alpha)}\|_{L_p}.
\end{equation}
The above inequality defines us the smallness of the initial data of $\nabla  \rho_i^0$, provided $\nabla \rho_i^{0(1-\alpha_p)} \in L_p$ and the $L_\infty$ norm of $\rho_i^0$ is sufficiently small.

We multiply the $m$-th equation of \eqref{subsys} by $a_m^{-\alpha}$ which gives
\begin{equation*}
\displaystyle  \frac{1}{1-\alpha} (a_m^{1-\alpha})_t + \frac{1}{1-\alpha} u \cdot \nabla (a_m^{1-\alpha}) =
  - \omega_m(\vec a)a_m^{-\alpha}.
\end{equation*}
Then we differentiate the above identity with respect to $x_j$, denoted again by $x$, getting 
\begin{equation} \label{3:6c}
\displaystyle  \frac{1}{1-\alpha} (a_m^{1-\alpha})_{xt} + \frac{1}{1-\alpha} u \cdot \nabla (a_m^{1-\alpha})_x+
\frac{1}{1-\alpha} u_x \cdot \nabla (a_m^{1-\alpha})
= - \de_x \Big(\omega_m(\vec a)a_m^{-\alpha}\Big).
\end{equation}
We aim at obtaining the bound in $L_p$, so we treat \eqref{3:6c} by $a_m^{-(p-1)\alpha} |a_{m,x}|^{p-2}a_{m,x}$. By the identity 
$$
\frac{1}{1-\alpha} a^{-\alpha(p-1)}|a_x|^{p-2}a_x (a^{1-\alpha})_{xt}
= \frac{1}{p(1-\alpha)^p} \de_t |(a^{1-\alpha})_x|^p
$$
we obtain
\begin{equation} \label{3:7}
 \frac{1}{p(1-\alpha)^p} \frac{d}{dt} \|(a_m^{1-\alpha})_x\|_{L_p}^p  \leq C\|\nabla u\|_{L_\infty} 
 \|\nabla (a_m^{1-\alpha})\|_{L_p}^p \\
 -\int_{\R^3} a_m^{-(p-1)\alpha} |a_{m,x}|^{p-2}a_{m,x}\de_x \Big(\omega_m(\vec a)a_m^{-\alpha}\Big)\,dx.
\end{equation} 
Now we sum the above inequalities over $m$ and variables $x_j$. By the second condition of \eqref{omega:pos2} the sum of the integral terms is again nonnegative, therefore 
\begin{equation*}
\frac{d}{dt} \Big( \sum_{m=1}^k \|\nabla a_m^{(1-\alpha)}\|_{L_p}^p \Big)
\leq C \|\nabla u\|_{\infty} \Big( \sum_{m=1}^k \|\nabla a_m^{(1-\alpha)}\|_{L_p}^p \Big)
\end{equation*}
with $C=C(p,\alpha)$, 
which implies 
\begin{equation} \label{3:14}
\sup_{t \in (0,T)}\left( \sum_{m=1}^k \|\nabla a_m^{(1-\alpha)}(t)\|_{L_p}^p \right)  
\leq \Big( \sum_{m=1}^k \|\nabla a_m^{(1-\alpha)}(0)\|_{L_p}^p \Big)
\exp\left(C\int_0^T \|\nabla u(s)\|_{\infty}\,ds\right).  
\end{equation}
Therefore by \eqref{est-rho}, \eqref{3:14}  and (\ref{yy1}) we obtain
\begin{equation}
\begin{aligned}
\sup_{t \in (0,T)}\|\nabla \vec a(t)\|_{L_p}^p & \leq C\|\vec a^{\;0}\|_{L_\infty}^{p\alpha}\Big( \sum_{m=1}^k \|\nabla a_m^{(1-\alpha)}(0)\|_{L_p}^p \Big)\exp\left(C\int_0^T \|\nabla u(s)\|_{\infty}\,ds\right)
\end{aligned}
\end{equation}
%
Having a bound on $\nabla \vec a$ we deal with $\nabla \vec b$.
We take advantage of the approach via the function $B$ introduced before. Since 
$$
\partial_t(B+\sum_{i=1}^k a_i) + 
u\cdot \nabla (B+\sum_{i=1}^k a_i)=0,
$$
by \eqref{B-zero} we deduce that\\
\begin{equation*}
\sup_t    \|\nabla B(t) + \nabla \sum_{i=1}^k a_k(t)\|_{L_p} \leq \|\nabla \vec a\,^0\|_{L_p}
\exp\left( \int_0^\infty \|\nabla u\|_{L_\infty} ds\right),
\end{equation*}
so the triangle inequality implies
\begin{equation*}
  \sup_t    \|\nabla B(t)\|_{L_p} \leq C\|\nabla \vec a\,^0\|_{L_p}
\exp\left( \int_0^\infty \|\nabla u\|_{L_\infty} ds\right).  
\end{equation*}
Therefore taking into account (\ref{b-0-B}) we have
\begin{equation} \label{3:15}
\sup_t \|\nabla b_j^{zero}(t)\|_{L_p} \leq C\theta_j\|\nabla \vec a\,^0\|_{L_p}\exp\left( \int_0^\infty \|\nabla u\|_{L_\infty} ds\right).    
\end{equation}
On the other hand, standard $L_p$ estimate for \eqref{b-ini} gives
\begin{equation} \label{3:16}
\sup_t \|\nabla b_j^{ini}(t)\|_{L_p} \leq C\|\nabla b_j^0\|_{L_p}\exp\left( \int_0^\infty \|\nabla u\|_{L_\infty} ds\right).    
\end{equation}
Combining \eqref{3:15} and \eqref{3:16} we get
\begin{equation*}
    \sup_t \|\nabla b_j^{ini}(t)\|_{L_p} +
    \sup_t \|\nabla b_j^{zero}(t)\|_{L_p} \leq C(\|\nabla b_j^0\|_{L_p}+\theta_j\|\nabla \vec a\,^0\|_{L_p}) \exp\left( \int_0^\infty \|\nabla u\|_{L_\infty} ds\right),
\end{equation*}
and so by (\ref{b-split}) we find
\begin{equation*}
    \sup_t \|\nabla \vec b(t)\|_{L_p} \leq C\|\nabla \vec b\,^0,\nabla \vec a\,^0\|_{L_p} \exp\left( \int_0^\infty \|\nabla u\|_{L_\infty} ds\right).
\end{equation*}
To close the part concerning the density estimates observe that Lemmas \ref{l:3} and \ref{l:4} together with velocity estimates imply 

\begin{equation} \label{pert:small}
\begin{aligned}
&\|\vec\rho - \hat e_1\|_{L_\infty(\R_+;L_\infty\cap \dot W^1_3 \cap \dot W^1_6)} \leq \\ 
&\phi\big( \|\vec a_0,\vec b_0\|_{\infty},\|\nabla a_0,\nabla b_0\|_{L_3\cap L_6}, \sum_{i=1}^k\|\nabla (a_i^0)^{(1-\alpha_3)}\|_{L_3}, \sum_{i=1}^k\|\nabla (a_i^0)^{(1-\alpha_6)}\|_{L_6} \big).
\end{aligned}
\end{equation}
\qed

{\bf Proof of the estimate \eqref{0:6}.}
We have now all necessary results to easily deduce the key a priori estimate \eqref{0:6}. Let us denote 
\begin{align*}
&A = \|u\|_{W^{2,1}_{2,(4/3,1)}\cap W^{2,1}_{5/4,(5/4,1)} }+\|tu\|_{W^{2,1}_{6,(4,1)}\cap  W^{2,1}_{2,(4,1)}
}+\|\nabla u\|_{L_1(\R_+;L_\infty)},\\
&D_0 = \left(\|\vec a_0,\vec b_0\|_{\infty},\|\nabla a_0,\nabla b_0\|_{L_3\cap L_6}, \sum_{i=1}^k\|\nabla (a_i^0)^{(1-\alpha_3)}\|_{L_3}, \sum_{i=1}^k\|\nabla (a_i^0)^{(1-\alpha_6)}\|_{L_6}\right).  
\end{align*}

Taking into account Lemmas \ref{l:2}, \ref{l:21}, \ref{lem:4} 
and \eqref{pert:small} we obtain 
$$
A \leq \phi(D_0)A^2 + \|u_0\|_{\dot B^{1/2}_{2,1}\cap \dot B^{2/5}_{5/4,1}},  
$$
which together with \eqref{pert:small} implies \eqref{0:6}.


\medskip 

%


{\bf Example.} Before we demonstrate an example of a simple reaction satisfying condition \eqref{omega:pos2} we shall make some comments about the latter. It seems quite technical, as we can see from the above proof we need it to assure proper sign of the right hand side. Possible generalizations of this condition are an interesting task for the future, the problem is related to the regularity of hyperbolic systems, therefore it is of independent interest. 
A stationary version of the below system with more general right hand side admitting reversibility of reactions, however with diffusion, was investigated in \cite{EZ}.

Let us now proceed with our example. Consider a simple reaction with two reactants and one product  
\begin{equation}
 a_1+a_2 \longmapsto b, \quad \omega_1(a_1,a_2)=a_1^2 a_2, \quad 
 \omega_2(a_1,a_2)=a_2^2 a_1.
\end{equation}
It is described by the following system:
\begin{equation} \label{toymodel}
 \begin{array}{l}
  a_{1,t} + u \cdot \nabla a_1 = -a_1^2 a_2,\\[5pt]
  a_{2,t} + u \cdot \nabla a_2 = - a_2^2 a_1, \\[5pt]
  b_t + u \cdot \nabla b = a_1^2 a_2 + a_2^2 a_1, \\[5pt]
  w_t + u\cdot \nabla w = 0.
 \end{array}
\end{equation}
We  show that above functions $\omega_1,\omega_2$ satisfy \eqref{omega:pos2}. 
%
We have 
$$
\omega_1(a_1,a_2)a_1^{-\alpha}=-a_1^{2-\alpha}a_2 \quad
\omega_2(a_1,a_2)a_2^{-\alpha}=-a_2^{2-\alpha}a_1,
$$
therefore 
\begin{align*}
& \de_x\big(\omega_1(a_1,a_2)a_1^{-\alpha}\big)=
-(2-\alpha)a_1^{1-\alpha}a_2a_{1,x}-a_1^{2-\alpha}a_{2,x},\\[5pt]
& \de_x\big(\omega_2(a_1,a_2)b^{-\alpha}\big)=
-(2-\alpha)a_2^{1-\alpha}a_1 a_{2,x}-a_2^{2-\alpha}a_{1,x},
\end{align*}
and so 
\begin{align*}
&a_1^{-(p-1)}\alpha|a_{1,x}|^{p-2}a_{1,x} \de_x\big(\omega_1(a_1,a_2)a_1^{-\alpha}\big) +
a_2^{-(p-1)}\alpha|a_{2,x}|^{p-2}a_{2,x}
\de_x\big(\omega_2(a_1,a_2)a_2^{-\alpha}\big) = \\[5pt]
&(2-\alpha)a_1^{1-p\alpha}a_2|a_{1,x}|^p+a_1^{2-p\alpha}a_{2,x}|a_{1,x}|^{p-2}a_{1,x}
+(2-\alpha)a_2^{1-p\alpha}a_1|a_{2,x}|^p+a_2^{2-p\alpha}a_{1,x}|a_{2,x}|^{p-2}a_{2,x}.
\end{align*}
To have a chance to show the nonnegativity of the above expression we put 
\begin{equation} \label{3:9}
2-p\alpha =0.
\end{equation}
Then the sum of the integrals read
\begin{equation} \label{3:10a}
\begin{aligned}
& \int (2-\alpha) a_1^{-1} a_2 |a_{1,x}|^p + a_{2,x} |a_{1,x}|^{p-2}a_{1,x} + a_{1,x} |a_{2,x}|^{p-2}a_{2,x} +
 (2-\alpha) a_2^{-1} a_1 |a_{2,x}|^p dx \\
& \geq |a_{2,x}|^p \int_{\{a_{2,x} \neq 0\}} \left\{ (2-\alpha) a_1^{-1} a_2 \left(\frac{|a_{1,x}|}{|a_{2,x}|}\right)^p + (2-\alpha) a_2^{-1}a_1 - \left(\frac{|a_{1,x}|}{|a_{2,x}|}\right)^{p-1} - \frac{|a_{1,x}|}{|a_{2,x}|}\right\}\,dx.     
\end{aligned}
\end{equation}
Denoting $\beta=a_1^{-1}a_2,\; \zeta = |a_{1,x}||a_{2,x}|^{-1}$ we see that the rhs of \eqref{3:10a} is nonnegative as 
\begin{equation} \label{3:11}
 (2-\alpha) \beta \zeta^p + (2-\alpha) \beta^{-1} - \zeta^{p-1} - \zeta \geq 0 \qquad {\rm for} \; \beta, \zeta > 0. 
\end{equation}
In order to show \eqref{3:11} we first observe that by Young inequality we have for any $\lambda>0$
\begin{equation} \label{3:12}
 \begin{array}{l}
\displaystyle \zeta = \zeta \lambda^{1/p'} \lambda^{-1/p'}  \leq \frac{1}{p} \lambda^{p-1} \zeta^p + \frac{1}{p'} \lambda^{-1},\\[5pt]
\displaystyle \zeta^{p-1} = \zeta^{p-1} \lambda^{1/p'} \lambda^{-1/p'}  \leq \frac{1}{p'} \lambda \zeta^p + \frac{1}{p} \lambda^{1-p}.
 \end{array}
\end{equation}
Taking $\lambda = \beta$ in \eqref{3:12} we obtain 
\begin{equation} \label{3:13}
 \zeta \leq \frac{1}{p} \beta^{p-1} \zeta^p + \frac{1}{p'} \beta^{-1},
\qquad
\zeta^{p-1} \leq \frac{1}{p'} \beta \zeta^p + \frac{1}{p} \beta^{1-p}.
\end{equation}
Applying \eqref{3:13} we find for $p\geq 2$ (which is equivalent to $\alpha\leq 1$ by \eqref{3:9}):
\begin{equation*}
\begin{aligned}
& (2-\alpha) \beta \zeta^p + (2-\alpha) \beta^{-1} - \zeta^{p-1} - \zeta \\ 
& \geq (2-\alpha)\beta \zeta^p + (2-\alpha)\beta^{-1} - \frac{1}{p} \beta^{p-1} \zeta^p - \frac{1}{p'} \beta^{-1}
 -\frac{1}{p'} \beta \zeta^p - \frac{1}{p} \beta^{1-p} \\
& \geq (2-\alpha)\beta \zeta^p + (2-\alpha)\beta^{-1} - \beta \zeta^p - \beta^{-1}\geq  0,
 \end{aligned}
\end{equation*}
which proves the second property from  \eqref{omega:pos2}.

\section{Existence} \label{sec:existence} 
Here we apply the estimates from the previous section to prove Theorem \ref{t1} applying the contraction principle. Note that the system has a hyperbolic character, hence a direct application of 
the Banach fixed point theorem does not work. We will adopt the technique from \cite{DMjfa}, in order to avoid the need of application of a tedious technique of the Lagrangian coordinates like in 
\cite{DM12,MX}.

Basing on the experience in the analysis of inhomogenous Navier-Stokes system \cite{DMjfa}, we construct the  approximate solutions by the following iteration
\begin{equation} \label{sysiter}
\begin{array}{lcl}
 \rho_{i,t}^{n+1} + u^n \cdot \nabla 
 \rho_i^{n+1} = \omega_i(\vec \rho\,^{n+1}),\quad i=1,\ldots,M & \mbox{ in } & \R^3 \times [0,T),\\[5pt]
 \rho^{n+1} u_t^{n+1} + \rho^{n+1} u^n\cdot \nabla u^{n+1} - \div \!\!\left(\nu (\vec  \rho\,^{n}) \D(u^{n+1})\right) + \nabla \pi^{n+1} =0& \mbox{ in } & \R^3 \times [0,T),\\[5pt]
 \div u^{n+1}=0 & \mbox{ in } & \R^3 \times [0,T).
\end{array}
\end{equation}
The initial data for ${\vec\rho}^{\;n+1}$ and $u^{n+1}$ are given by (\ref{init}).
The idea is the following, having $u^n$ we solve the first equations and get $\vec\rho^{\;n+1}$, then we solve the momentum equation and get $u^{n+1}$. Without loss of generality we put $u^0\equiv 0$.

The solvability of the above system is a classical result, as it can be treated as a simple perturbation of a linear system: nonlinearity in \eqref{sysiter}$_1$ is in a regular lower order term $\omega_i$, since $u^n$ is given we just solve the transport equation. While \eqref{sysiter}$_2$ is linear in $u^{n+1}$, it is a perturbation of the classical Navier-Stokes system. The a priori estimate yields the scheme of proving existence using the Banach principle.

Moreover we observe that the estimates from section 4  
imply the following a priori estimate for $(u^n,\rho^n)$:
\begin{multline}\label{est-app}
  \sup_{n \in \mathbb{N}} \left(  \| \vec \rho\,^n -\hat{e}_1\|_{L_\infty(\R_+;L_\infty \cap\dot W^1_3\cap \dot W^1_6)}+\|u^n\|_{W^{2,1}_{2,(4/3,1)}(\R^3 \times \R_+)} +
  \|u^n\|_{W^{2,1}_{5/4,(5/4,1)}(\R^3 \times \R_+)}\right.\\
  \left.
  +\|tu^n\|_{W^{2,1}_{2,(4,1)}(\R^3 \times \R_+)}
  +\|tu^n\|_{W^{2,1}_{6,(4,1)}(\R^3 \times \R_+)}
  \right) \leq \epsilon.
\end{multline}
Because of the hyperbolic character of the transport equation we are not able to obtain the convergence of the approximative sequence in the topology determined by estimate (\ref{est-app}), therefore a direct application of the Banach fixed point theorem is impossible. Hence, following the typical trick for compressible or inhomogeneous Navier-Stokes equations, we consider the differences
\begin{equation}
    \delta \vec \rho\,^{n+1}:= \vec \rho\,^{n+1} - \vec \rho\,^n, 
    \qquad \delta u^{n+1}:= u^{n+1} - u^n.
\end{equation}
They satisfy the following system
\begin{equation} \label{Banach}
\begin{array}{lcl}
 \delta \rho_{i,t}^{n+1} + u^n \cdot \nabla 
 \delta \rho_i^{n+1} = \omega_i(\vec \rho\,^{n+1})
 - \omega_i(\vec \rho\,^n) - \delta u^n \cdot \nabla \rho^n_i,\quad i=1,\ldots,M & \mbox{ in } & \R^3 \times [0,T),\\[5pt]
 \rho^{n+1} \delta u_t^{n+1} + \rho^{n+1} u^n\cdot \nabla \delta u^{n+1} - \div \!\!\left(\nu (\vec  \rho\,^{n+1}) \D(\delta u^{n+1})\right) + \nabla \delta p^{n+1} \qquad  \,= & & \\[5pt]
\qquad \, -\delta \rho^{n+1} u_t^n 
-(\rho^{n+1} u^n - \rho^n u^{n-1})\cdot \nabla u^n 
+\div( (\nu(\vec \rho\,^{n+1}) - \nu( \vec \rho\,^n)) \D(u^n) ) & \mbox{ in } & \R^3 \times [0,T),\\[5pt]
 \div \delta u^{n+1}=0 & \mbox{ in } & \R^3 \times [0,T).
\end{array}
\end{equation}
We aim at showing the contraction in a larger space, namely we show
\begin{multline}\label{ee1}
    \|t^{-1/2}\delta \vec \rho\,^{n+1}\|_{L_\infty(0,T;L_2)}+
    \|\delta u^{n+1}\|_{L_\infty(0,T;L_2)}+
    \|\nabla \delta u^{n+1}\|_{L_2(0,T;L_2)} 
    \\
    \leq 
    \frac12  \left(\|t^{-1/2}\delta \vec \rho\,^{n}\|_{L_\infty(0,T;L_2)}+
    \|\delta u^{n}\|_{L_\infty(0,T;L_2)}+
    \|\nabla \delta u^{n}\|_{L_2(0,T;L_2)}
    \right). 
\end{multline}
Testing $(\ref{Banach})_1$ by $\delta \rho_{i}^{n+1}$ we get
\begin{equation*}
\begin{aligned}
    \frac12 \frac{d}{dt} \|\delta \vec \rho\,^{n+1}\|^2_{L_2} &\leq C_0\|\delta \vec \rho\,^{n+1}\|^2_{L_2}+ \|\delta u^n\|_{L^6} \|
    \nabla \vec \rho\,^n\|_{L_3}\|\delta \vec \rho\,^{n+1}\|_{L_2}\\[5pt]
    &\leq C_0\|\delta \vec \rho\,^{n+1}\|^2_{L_2}+ \ep\|\delta u^n\|_{L^6} \|\delta \vec \rho\,^{n+1}\|_{L_2}.
\end{aligned}
\end{equation*}
The first term of the RHS comes from the difference of $\omega$'s, from that reason we can estimate it with a uniform in $n$ and time constant $C_0$, related to the Lipschitz constant of $\vec \omega(\cdot)$. Next, we divide by $\|\delta \vec \rho\,^{n+1}\|_{L_2}$ and apply the Gronwall inequality
getting
\begin{equation*} \label{ee1a}
    \|\delta \vec \rho\,^{n+1}\|_{L_2}(t) \leq C \epsilon e^{C_0t} \int_0^t
    e^{-C_0 s} \|\delta u^n\|_{L_6}(s) ds \\
    \leq C \epsilon  e^{C_0T} t^{1/2} (\int_0^t \|\delta u^n\|_{L_6}^2(s)ds)^{1/2},
\end{equation*}
which implies for finite $T$
\begin{equation}\label{ee22}
    \sup_{t\in [0,T]} t^{-1/2}\|\delta \vec \rho\,^{n+1}\|_{L_2}(t)\leq C\epsilon e^{C_0T}\|\nabla \delta u^n\|_{L_2(0,T;L_2)}.
\end{equation}
Note that since $n$ is given in (\ref{ee22}) one can consider the limit $T\to 0^+$, 
showing that the limit is well posed and it is equal zero since $\|\nabla u^n\|_{L_2(0,T;L_2)}$ goes to zero as $T\to 0^+$.

In the same manner we consider the  momentum equation. Testing it by $\delta u^{n+1}$ we get the bound on $\delta u^{n+1}$. On the RHS we find four terms. One of them is
\begin{multline*}
     \left| \int_{\R^3} \delta \rho^{n+1}\, u^n_t\, \delta u^{n+1} dx\right|=
    \left|\int_{\R^3} t^{-1/2}\delta \rho^{n+1} \, t^{1/2} u^n_t \, \delta u^{n+1} dx\right| 
    \\
    \leq
    \|t^{-1/2}\delta \rho^{n+1}\|_{L_2}
    \|t^{1/2} u^n_t\|_{L^3}\|\delta u^{n+1}\|_{L_6}. \qquad
\end{multline*}
Now, let us note that $u_t \in L_{4/3}(0,T;L_2)$
and $tu_t \in L_4(0,T;L_6)$, so 
\begin{equation*}
    \|t^{1/2}u_t\|_{L_3} \leq \|u_t\|_{L_2}^{1/2}
    \|tu_t\|_{L_6}^{1/2}.
\end{equation*}
Thus by (\ref{est-app}) we find
\begin{equation*}
    \left\| \|t^{1/2}u_t\|_{L_3}(t) \right\|_{L_2} \leq 
    \|u_t\|_{L_{4/3}(0,T;L_2)}^{1/2} \|t u_t\|_{L_{4}(0,T;L_6)}^{1/2} \leq \epsilon.
\end{equation*}
The above estimate leads to the following bound
\begin{multline*}
     \left| \int_0^T\int_{\R^3} \delta \rho^{n+1}\, u^n_t \, \delta u^{n+1} dxdt \right| 
    \\[8pt]
    \leq
    \sup_t \|t^{-1/2}\delta \rho^{n+1}\|_{L_2}
    \|t^{1/2} u^n_t\|_{L_2(0,T;L_3)}\|\delta u^{n+1}\|_{L_2(0,T;L_6)} 
    \\
    \leq \epsilon \sup_t \|t^{-1/2}\delta \rho^{n+1}\|_{L_2}
    \|\delta u^{n+1}\|_{L_2(0,T;L_6)} .
\end{multline*}
Next we have a term which is easily bounded by Lemma \ref{lem:4}
\begin{multline*}
    \left| \int_0^T \int_{\R^3} \rho^n \,\delta u^n \cdot \nabla u^{n} \, \delta u^{n+1} dx dt \right| \leq 
    \\
    C \sup_t \|\rho^n\|_{L_\infty} \|\delta u^{n}\|_{L_\infty(0,T;L_2)}
    \|\delta u^{n+1}\|_{L_\infty(0,T;L_2)} \int_0^\infty 
    \|\nabla u^n\|_{L_\infty} dt \\
    \leq C \epsilon \|\delta u^{n}\|_{L_\infty(0,T;L_2)}
    \|\delta u^{n+1}\|_{L_\infty(0,T;L_2)}. 
\end{multline*}
Next
\begin{equation*}
    \left|\int_0^T \int_{\R^3} \delta \rho^{n+1} \, u^n \cdot \nabla u^n \, \delta u^{n+1} dxdt\right| \\
    \leq 
    \epsilon \sup_t \|t^{-1/2}\delta \rho^{n+1}\|_{L_2}
    \|\delta u^{n+1} \|_{L_2(0,T;L_6)},
\end{equation*}
since 
\begin{equation*}
    \|t^{1/2} u^n \cdot \nabla u^n\|_{L_2(0,T;L_3)}\leq C\epsilon^2.
\end{equation*}
In order to obtain the above estimate note that by the standard embedding $\nabla W^{2,1}_{4,(4/3,1)} \subset L_2(\R_+;L_3)$ and by interpolation we get
\begin{equation*}
\|t^{1/2}u\|_{L_\infty} \leq C\|u\|_{\dot B^{-1/2}_{6,1}}^{1/2} \|tu\|_{\dot B^{3/2}_{6,1}}^{1/2}
\leq C\epsilon. 
\end{equation*}
With the last term with $\D(u^n)$ we proceed in the same manner.
Hence finally we deduce
\begin{multline} \label{ee1b}
    \|\delta u^{n+1}\|_{L_\infty(0,T;L_2)}^2+
    \|\nabla \delta u^{n+1}\|_{L_2(0,T;L_2)}^2 
    \\
    \leq 
    C\epsilon  \left(\|t^{-1/2}\delta \vec \rho\,^{n+1}\|_{L_\infty(0,T;L_2)}^2+
    \|\delta u^{n}\|_{L_\infty(0,T;L_2)}^2+
    \|\nabla \delta u^{n}\|_{L_2(0,T;L_2)}^2
    \right).  
\end{multline}
Taking $\epsilon$ small enough and $T$ fixed, in this framework we need to assume that $\epsilon$ is small even for small $T$, by (\ref{ee22}) we find the desired contraction condition (\ref{ee1}).
But $T$ is fixed from now. 
%
%
We  now proceed with convergence of the sequence of approximations.
Estimate \eqref{ee1} implies that as $n \to \infty$ (up to a sub-sequence) 
\begin{equation}\label{ee2}
    \vec \rho\,^n \to \vec \rho \mbox{ \ \ in \ } L_\infty(0,T;L_2) 
    \mbox{ \  \  and \ \ } 
    u^n \to u \mbox{ \ \ in \ }
    L_\infty(0,T;L_2) \cap L_2(0,T;\dot H^1)
\end{equation}
for some $(\vec \rho, u)$ in the spaces defined by the convergence. 

What we may say about the limit? The convergence is strong in the norm, so up to a sub-sequence we have the point-wise convergence. From the uniform estimate (\ref{est-app}), one could get also information in higher regularity
for the weak limits of the sequences. But the problem is that the Lorentz space $L_{p,1}(0,T,X)$ is not reflexive, so we can not find a weak limit in this regularity. For this reason we proceed as follows:

Note that from (\ref{est-app}) and (\ref{ee2}) we conclude that
\begin{equation}
    \vec \rho\,^n \rightharpoonup \vec \rho \mbox{ weakly} \ast \mbox{  in } L_\infty(0,T;L_\infty) \mbox{ \  and  \ } \vec \rho\,^n \to  \vec \rho \mbox{ strongly in } L_p(0,T;L_p) \mbox{ for any } p<\infty.
\end{equation}
It suffices to pass to the limit (for  a suitable subsequence) 
in (\ref{sysiter}) and get
\begin{equation} \label{sys4}
\begin{array}{lcl}
 \rho_{i,t} + u \cdot \nabla \rho_i = \omega_i(\vec \rho\,),\quad i=1,\ldots,M & \mbox{ in } & \R^3 \times [0,T),\\[5pt]
 \rho u_t  - \div \!\!\left(\nu (\vec  \rho\,) \D(u)\right) + \nabla \pi =-\rho u\cdot \nabla u & \mbox{ in } & \R^3 \times [0,T),\\[5pt]
 \div u=0 & \mbox{ in } & \R^3 \times [0,T).
\end{array}
\end{equation}

The key problem is to recover the regularity of the velocity in Lorentz spaces. The main term is the nonlinear one $\rho u\cdot \nabla u $. We put it on the RHS. There is a need to show that it belongs to 
$L_{4/3,1}(0,T;L_2)$. Since $\rho$ is bounded we just put our attention on $u\cdot \nabla u$. 
Based on estimates (\ref{est-app}) we conclude that
\begin{equation*}
    \|u^n\|_{L_{4,1}(0,T;L_6)} + \|\nabla u^n\|_{L_{2,1}(0,T;L_3)} \leq C \epsilon.
\end{equation*}
But by imbeddings \eqref{imbed} we have also
\begin{equation*}
    \|u^n\|_{L_{4,2}(0,T;L_6)} + \|\nabla u^n\|_{L_{2,2}(0,T;L_3)} \leq C \epsilon,
\end{equation*}
and now the spaces $L_{4,2}(0,T;L_6)$ and $L_{2,2}(0,T;L_3)$
are reflexive \cite{Tr}, so the weak limit of the velocity is an element of these spaces, i.e.
\begin{equation*}
    u \in L_{4,2}(0,T;L_6) \mbox{ \ \ and  \ \ } \nabla u \in L_{2,2}(0,T;L_3).
\end{equation*}
Now we apply a magic property of Lorentz spaces that in some sense nonlinearity improves the regularity, of course, only for the second index of the space.
Hence from the H\"older inequality we find that
\begin{equation*}
    \|\rho u \cdot \nabla u\|_{L_{4/3,1}(0,T;L_2)}
    \leq \|\rho\|_{L_\infty(0,T;L_\infty)}
    \|u\|_{L_{4,2}(0,T;L_6)}\|\nabla u\|_{L_{2,2}(0,T;L_3)}.
\end{equation*}
So we  restate the system (\ref{sys4}) in such a way that
\begin{equation}\label{sys5}
    \rho u_t  - \div \!\!\left(\nu (\vec  \rho\,) \D(u)\right) + \nabla \pi \in L_{4/3,1}(0,T;L_2).
\end{equation}
Taking into account the features of $\vec \rho$ we treat (\ref{sys5}) as a perturbation of the Stokes system with the RHS in $L_{4/3,1}(0,T;L_2)$, so we obtain the solvability in this framework with the bound
\begin{equation*}
    \|u\|_{W^{2,1}_{2,(4/3,1)} }\leq C\|\rho u \cdot \nabla u\|_{L_{4/3,1}(0,T;L_2)}+\|u_0\|_{\dot B^{1/2}_{2,1}}.
\end{equation*}
The estimates for $u$ in $W^{2,1}_{5/4,(5/4,1)}$, $tu$ in $W^{2,1}_{2,(4,1)}$ and  $tu$ in $W^{2,1}_{6,(4,1)}$ are obtained in the same way. It gives the existence on the time interval $[0,T]$. Repeating this procedure we restart the solution from $t=T$ and show the existence over time interval $[T,2T]$ and so on. The a priori estimates are valid for all times, in particular we control the norms of required time traces
allowing to initiate the next part of the solution for the next time step.
As the bounds are time independent  one can extend the solution to the whole halfline. The uniqueness follows from the method of the proof of the existence. Theorem \ref{t1} is therefore proved.
\qed

\noindent
{\bf Acknowledgements.} The authors have been partially supported by National Science Centre grant
No2018/29/B/ST1/00339 (Opus). The authors would also like to thank the anonymous Referees who pointed out several important points in the first version of the manuscript. Their careful lecture of the manuscript contributed a lot to the quality of the paper. 

\section*{Declarations}

{\bf Conflict of interest:} the authors declare no conflict of interest.

\vskip5mm

\noindent
{\bf Availability of data and materials:} does not apply to this work.

\end{document}